\documentclass[12pt]{amsart}

\usepackage{amsfonts,amsmath,amssymb}

\usepackage{enumerate}

\numberwithin{equation}{section}

\newtheorem{theorem}{Theorem}[section]

\newtheorem{lemma}[theorem]{Lemma}

\newtheorem{example}[theorem]{Example}

\newtheorem{definition}{Definition}[section]

\newtheorem{corollary}[theorem]{Corollary}

\newtheorem{remark}[theorem]{Remark}

\newcommand{\cl}[1]{\mathcal{#1}} 

\newcommand{\nor}[1]{\left\Vert #1\right\Vert}

\newcommand{\Lat}[1]{\mathrm{Lat}(#1)} 

\newcommand{\Alg}[1]{\mathrm{Alg}(#1)} 

\newcommand{\wsp}[1]{\overline{\mathrm{span}}^{\text{w*}}(#1)} 

\newcommand{\lsp}[1]{\overline{\mathrm{span}}(#1)} 

\begin{document}

\title{Morita equivalence of nest algebras}

\author[G.~K.~Eleftherakis]{G. K. Eleftherakis}

\address{Department of Mathematics, University of Athens,
Panepistimioupolis 157 84, Athens, Greece.}

\email{gelefth@math.uoa.gr}



\date{}

\maketitle

\begin{abstract}
Let $\cl N_1$ (resp. $\cl N_2$) be a nest, $A$ (resp. $B$) be the corresponding 
nest algebra, $A_0$ (resp. $B_0$) be the subalgebra of compact operators. We 
prove that the nests $\cl N_1, \cl N_2$ are isomorphic if and only if the algebras 
$A, B$ are weakly$-*$ Morita equivalent if and only if the algebras $A_0, B_0$ are strongly 
Morita equivalent. We characterize the nest isomorphisms which implement stable 
isomorphism between the corresponding nest algebras.
\end{abstract}

\section{Introduction}

Rieffel, introduced the idea of Morita equivalence in Operator Theory developing 
the theory of Morita equivalence for $C^*$ and $W^*$ algebras \cite{rief}. 
After the advent of the theory of operator spaces and operator algebras 
a parallel Morita theory for non-selfadjoint algebras was developed by 
Blecher, Muhly and Paulsen \cite{bmp}, \cite{blecher}. We call this equivalence strong Morita 
equivalence. 

Recently, two approaches have been suggested for the Morita equivalence of dual operator algebras.
 The one was introduced \cite{ele1}  by the author of this article and 
it is equivalent to the notion of stable isomorphism of dual operator algebras \cite{elepaul}. 
We call this equivalence $\Delta -$equivalence. The other was introduced by Blecher and 
Kashyap \cite{bk}, \cite{kashyap} and it is strictly weaker than $\Delta -$equivalence. 
This equivalence is called weak$-*$ Morita equivalence. It is interesting that if $A$ and 
$B$ are strongly Morita equivalent approximately unital operator algebras 
then the second dual operator algebras $A^{**}, B^{**}$ are weakly$-*$ Morita equivalent 
\cite{bk}. New results on weak$-*$ Morita equivalence and $\Delta -$equivalence can be found 
in \cite{bkraus}.

In this paper we prove that strong and weak$-*$ Morita equivalence is a lattice property 
for nest algebras. Particularly we prove that if $A$ and $B$ are nest algebras 
and $A_0, B_0$ are the subalgebras of compact operators then $A_0$ and $B_0$ are 
strongly Morita equivalent if and only if  $A$ and $B$ are weakly$-*$ Morita equivalent 
if and only if the nests $ \mathrm{Lat}(A), \mathrm{Lat}(B) $  are isomorphic. The main tool of the proof 
is that if $\theta : \mathrm{Lat}(A)\rightarrow  \mathrm{Lat}(B) $ is a nest isomorphism 
we can construct a dual operator $A-B$ bimodule $Y$ and a dual operator $B-A$ bimodule $X$ 
such that the identity operator of $A$ is the limit in strong operator topology 
of a net of finite rank contractions $(f_\lambda) $ where 
every $f_\lambda $ is the norm limit of a sequence $(y_i^\lambda x_i^\lambda )_{i\in 
\mathbb{N}},$  where $y_i^\lambda $ is a contractive row operator with finite entries from $Y$  and 
  $x_i^\lambda $ is a contractive column operator with finite entries from $X.$ Similarly we 
can decompose the identity of the algebra  $B.$ This can be considered as a generalization 
of the Erdos density Theorem for nest algebras \cite{dav}.

In section 3 we prove that two nest algebras are weakly$-*$ Morita equivalent if and only if 
they are 
spatially Morita equivalent (definition \ref{spat1}). Also we prove that every spatially Morita equivalent dual operator algebra 
with a nest algebra is weakly$-*$ Morita equivalent with this nest algebra. It is interesting 
that this does not happen  
for the more general class of operator algebras, the CSL algebras.

In section 4 we present a measure-theoretic result which describes when two 
separably acting nest algebras are stably isomorphic. As it was pointed out in \cite{bk} 
the \cite[example 3.7]{ele2} is an example of weak$-*$ Morita equivalent algebras 
which are not stably isomorphic. Using the results of this paper 
we give a new proof of the fact that weak$-*$ Morita equivalence is strictly weaker 
than $\Delta -$equivalence. 

In section 5 we present a counterexample which states that the second duals of 
two unital strongly Morita equivalent algebras are not necessarily stably isomorphic. 

\bigskip

In what follows we describe the notions we use in this paper. Since we use extensively 
the basics of Operator Space Theory, we refer the reader to the monographs \cite{bm}, 
\cite{er}, \cite{paul} and \cite{pisier}  
for further details. A (normal) \textbf{representation} 
of a (dual) operator algebra $A$ is a ($w^*-$continuous) completely contractive homomorphism 
 $\alpha : A\rightarrow B(H)$ on a Hilbert space $H.$ In the case $A$ is unital, we assume that 
$\alpha $ is unital.

Let $H, K$ be Hilbert spaces and $A\subset B(H)$ be an algebra. A subspace $X\subset B(K, H)$ 
is called a left module over $A$ if $AX\subset X.$ Similarly we can define 
the right modules over $A.$ A left and right 
module over $A$ is called a bimodule over $A.$ 
An abstract left (right) operator module over an abstract operator 
algebra $A$ is an operator space $Y$ such that there exist a 
 completely contractive bilinear map $A\times Y\rightarrow Y \;\;(Y\times A\rightarrow Y).$
 A left and right operator 
module over $A$ is called an operator bimodule over $A.$ 
 
 If $A$ is a dual 
operator algebra and $Y$ is a dual operator space we say that $Y$ is a left (right) 
dual operator module 
if the above completely contractive bilinear map is separately $w^*-$continuous. A left and right dual operator 
module over $A$ is called a dual operator bimodule over $A.$ 

Two operator bimodules $Y$ and $Z$ over an operator algebra $A$ are 
called isomorphic as operator bimodules if there exists a completely 
isometric and onto $A-$module map $\pi : Y\rightarrow Z.$ We 
denote $Y\cong Z$ as operator bimodules. In the case $A$ is a dual operator algebra 
and $Y, Z$ are dual operator bimodules we denote $Y\cong Z$ as dual operator bimodules if the 
above completely isometric and onto $A-$module map  
$\pi $ is $w^*-$(bi)continuous.

If $Y$ is a right operator module over an operator algebra $A$ and $X$ 
is a left operator module over $A$ we denote by $Y\otimes ^h_A X$ the balanced 
Haagerup tensor product of $Y$ and $X$  which linearizes the completely bounded $A-$balanced 
bilinear maps \cite[3.4]{bm}. If $Y$ (resp. $X$) is a left (resp. right) operator 
module over an operator algebra $B$ then $Y\otimes ^h_A X$ is also a left (resp. right) 
operator module over $B,$ \cite[Lemma 2.4]{bmp}.

If $Y$ is a dual right operator module over a dual operator algebra $A$ and $X$ is 
a left dual operator module over $A$ we denote by $Y\otimes ^{\sigma h}_A X$ the balanced 
normal Haagerup tensor product of $Y$ and $X$ which linearizes the separately 
$w^*-$continuous completely bounded $A-$balanced 
bilinear maps \cite{elepaul}. In the case $Y$ (resp. $X$) is a left (resp. right) dual 
operator module 
over a dual operator algebra $B$ then $Y\otimes ^{\sigma h}_A X$ is also a left 
(resp. right) dual operator module over $B,$ \cite{elepaul}.

We give now the two definitions of Morita equivalence using in this paper:

\begin{definition}\label{strong}\cite{bmp} The operator algebras $A, B$ are called 
\textbf{strongly Morita equivalent} if there exist an $A-B$ operator module $X$ and 
a $B-A$ operator module $Y$ such that $A\cong X\otimes ^h_BY$ and $B\cong Y\otimes ^h_AX$ 
as $A$ and $B$ operator bimodules respectively. 
\end{definition}

\begin{definition}\label{weak}\cite{bk} The dual operator algebras $A, B$ are called 
\textbf{weakly$-*$ Morita equivalent} if there exist an $A-B$ dual operator module $X$ and 
a $B-A$ dual operator module $Y$ such that $A\cong X\otimes ^{\sigma h}_BY$ and 
$B\cong Y\otimes ^{\sigma h}_AX$  as $A$ and $B$ dual operator bimodules respectively. 
\end{definition}

If $X$ is a subspace of $B(H, K),$ where $H$ and $K$ are Hilbert spaces, we denote by 
$R^{fin}_\infty (X)$ (resp. $C^{fin}_\infty (X)$) 
the space of operators $(x_1, x_2,...): H^\infty \rightarrow K$ (resp. $(x_1, x_2,...)^T: 
H \rightarrow K^\infty $) such that $x_i\in X$ for all $i$ and there exists $n_0\in 
\mathbb{N}$ such that $x_n=0$ for all $n\geq n_0.$ 

If $s_1= (s_1^1,s^1_2,...,s^1_{n_1},0,0,...), s^1_{n_1}\neq 0 $ and 
$s_2=( s_1^2,s^2_2,...,s^2_{n_2},0,0,...), s^2_{n_2}\neq 0 $ 
are operators in $R^{fin}_\infty (X)$ we denote by $(s_1, s_2)$ the operator 
$$ (s_1^1,s^1_2,...,s^1_{n_1},s_1^2,s^2_2,...,s^2_{n_2},0,0,...) $$ 
which also belongs to $R^{fin}_\infty (X).$ In the same way if $s_1, s_2,...,s_n\in 
R^{fin}_\infty (X) $ we define the operator $(s_1, s_2,...,s_n)\in 
R^{fin}_\infty (X). $ Similarly if $t_1, t_2,...,t_n\in 
C^{fin}_\infty (X) $ we define the operator $(t_1, t_2,...,t_n)^T\in 
C^{fin}_\infty (X). $ 

A \textbf{nest} $\cl N$ is a totally ordered set of projections of a Hilbert space 
$H$ containing the zero and identity operators which is closed under arbitrary 
intersections and closed spans. The corresponding \textbf{nest algebra} is 
$$ \mathrm{Alg}(\cl N) =\{x\in B(H): N^\bot xN=0\;\;\forall \;\;N\;\;\in \cl N\}.$$
If $N\in \cl N$ we denote by $N_-$ the projection onto the closed span of the 
union $\cup _{\stackrel{M<N}{M\in \cl N}}(M(H)).$ If $N_-<N$ we call the projection 
$N\ominus N_-$ an \textbf{atom}. 
If the nest has not atoms is called a \textbf{continuous nest}. If the atoms span 
the identity operator the nest is called a \textbf{totally atomic nest}. An order preserving 
1-1 and onto map  
between two nests is called a \textbf{nest isomorphism}.

If $\cl N_1$ and $\cl N_2$ are nests acting on the Hilbert spaces $H_1, H_2$ 
respectively and $\theta : \cl N_1\rightarrow \cl N_2$ is a nest isomorphism  
we denote by $Op(\theta )$ the space of operators $x\in B(H_1, H_2)$ 
satisfying $\theta (N)^\bot x N=0$ for all $N\in \cl N_1.$ Observe that 
$Op(\theta )$  is an $ \mathrm{Alg}(\cl N_2) - \mathrm{Alg}(\cl N_1) $ bimodule.

Finally, if $X$ is a normed space we denote by $Ball(X)$ the unit ball of $X$ 
and by $X^*$ its dual space.  
If $H_1, H_2$ are Hilbert spaces and $\xi \in H_2, \eta \in H_1$ are vectors 
we denote by $\xi \otimes \eta^* $ the rank 1 operator sending every 
$\omega \in H_1$ to $<\omega ,\eta >\xi \in H_2,$ where $<\cdot ,\cdot>$ is 
the inner product of $H_1.$  
 Also we symbolize the  strong operator topology by SOT.

\section{Morita equivalence for nest algebras}

In this section we fix nests $\cl N_1, \cl N_2$ acting on the Hilbert 
spaces $H_1, H_2$ respectively, and a nest isomorphism $\theta : \cl N_1\rightarrow 
\cl N_2.$ We denote $ A=\Alg{ \cl N_1}, B=\Alg{ \cl N_2}, X=\mathrm{Op}(\theta ), 
Y=\mathrm{Op}(\theta^{-1} ).$ If $Z$ is a space of operators we denote its subspace of 
compact operators by $Z_0.$ Observe that 
$$AYB\subset Y,\;\; BXA\subset X,\;\; YX\subset A, \;\;XY\subset B,$$
$$A_0Y_0B_0\subset Y_0,\;\; B_0X_0A_0\subset X_0,\;\; Y_0X_0\subset A_0,\;\; X_0Y_0\subset B_0.$$

The main result of this section is Theorem \ref{2.11}. In particular we are 
going to prove that 
$$ A_0\cong Y_0\otimes ^h_{B_0}X_0,\;\;  B_0\cong X_0\otimes ^h_{A_0}Y_0,\;\; 
A\cong Y\otimes ^{\sigma h}_BX,\;\;  B\cong X\otimes ^{\sigma h}_AY.$$
Suppose that $p=\vee \{N\ominus N_-: N\in \cl N_1\}.$ The following lemmas are 
used in  Theorem \ref{2.6}, where we are going to prove a variant of the Erdos density Theorem for 
nest algebras: 
There exists a net of finite rank contractions $(f_\lambda )\subset A$ converging 
in SOT topology to the identity operator of $H_1,$ where every $f_\lambda $ is the norm 
limit of a sequence $(y_i^\lambda x_i^\lambda )_{i\in 
\mathbb{N}}$ where $ y_i^\lambda \in Ball(R_\infty ^{fin}(Y_0)), x_i^\lambda \in Ball(C_\infty 
^{fin}(X_0))$ for all $i,\lambda .$ 

\begin{lemma}\label{2.1} There exists a net $(l_\lambda )$ of finite rank contractions converging 
 in SOT topology to the projection $p$ such that $l_\lambda =s_\lambda t_\lambda $ 
where $ s_\lambda \in Ball(R_\infty ^{fin}(Y_0)),$ $ t_\lambda \in Ball(C_\infty 
^{fin}(X_0))$ for all $\lambda .$
\end{lemma}
\textbf{Proof} Suppose that $p=\vee _{k\in J}p_k$ where $p_k=N_k\ominus (N_k)_-$ 
for $N_k\in \cl N_1, k\in J.$ Choose a net of finite rank contractions $(f_i)_{i\in I}$ 
converging in SOT topology to the identity operator of $H_1.$ If $\cl F=\{F: F\text{\;\;finite
\;\;subset\;\;of\;\;}J\}$ the family $(g_{F,i})_{(F,i)}$ indexed by $\cl F\times I$ 
where $g_{F,i}=\sum_{k\in F}p_kf_ip_k$ is a net. Observe that every  $g_{F,i}$ is 
a finite rank contraction belonging to $A.$ We can easily check that 
$SOT-\lim_{(F,i)}g_{F,i}=\vee _kp_k=p.$

Let $f=p_kf_ip_k$ for some $k\in J, i\in I$ with polar decomposition $f=u|f|.$ Suppose that 
$$|f|=\sum_{j=1}^n\lambda _j
\xi_j\otimes  \xi_j^* $$ for $\lambda _j\geq 0$ and $\xi _j$ orthogonal vectors of $p_k(H_1).$
Choose a unit vector $\eta$  in $(\theta (N_k)\ominus \theta (N_k)_-)(H_2).$ Now we have 
$$|f|=\sum_{j=1}^n\lambda _j\xi_j\otimes \eta ^*\cdot\eta \otimes  \xi_j^* =yy^*$$ 
where $y=( \sqrt{\lambda _1}\xi _1\otimes \eta^* ,...,\sqrt{\lambda _n}\xi _n\otimes \eta^* ).$
Observe that $f=uyy^*$ and $uy\in Ball(R_\infty ^{fin}(Y_0)), y^*\in 
Ball(C_\infty ^{fin}(X_0)).$ 

Suppose now that $F=\{j_1,...,j_n\}\subset J, i\in I$ and 
$g_{F,i}=\sum_{k=1}^n p_{j_k}f_ip_{j_k} .$ By the above arguments $p_{j_k}f_ip_{j_k}=
s_kt_k$ where $s_k \in Ball( R_\infty ^{fin}(Y_0) ), t_k \in 
Ball( C_\infty ^{fin}(X_0) ).$ So $g_{F,i}=st$ where 
$$s=(s_1,...,s_n)\in R_\infty ^{fin}(Y_0) ,\;\; t=(t_1,...,t_n)^T\in C_\infty ^{fin}(X_0).$$ 
Also, since the projections $(p_k)_{k\in J}$ are pairwise orthogonal and $\|s_ks_k^*\|\leq 1$ for 
all $k$ we have that $$\|s\|^2= \|\sum_{k=1}^n s_ks_k^*\| 
=\|\sum_{k=1}^np_{j_k}s_ks_k^*p_{j_k}\|\leq 1 .$$
 Similarly we can prove $\|t\|\leq 1$ and this completes the proof.$\qquad \Box$

\begin{lemma}\label{2.2} Suppose that $p^\bot \neq 0$ $\xi, \eta \in Ball(H_1)$ and 
$N\in \cl N_1$ such that $\xi =p^\bot N(\xi ), \eta =p^\bot N^\bot (\eta ).$ 
There exist rank $1$ operators $ (s_n)_{n\in \mathbb{N}}\subset Ball(Y), 
(t_n)_{n\in \mathbb{N}}\subset Ball(X), $ such that the operator $\xi\otimes  \eta^* $ 
is the norm limit of the sequence $(s_nt_n)_{n\in \mathbb{N}}.$ 
\end{lemma}
\textbf{Proof} We define the continuous order preserving map 
$$\phi : p^\bot \cl N_1 \rightarrow [0,\|\xi \|^2]: p^\bot M\rightarrow 
\|p^\bot M(\xi )\|^2.$$ The nest $p^\bot \cl N_1$ is continuous, so $\phi $ 
is onto $[0, \|\xi \|^2 ].$ 
Choose a strictly increasing sequence $(\lambda _n)$ such that $\lambda _n\rightarrow \|\xi \|^2 
.$ Choose $N_n\in \cl N_1$ such that $\phi ( p^\bot N_n )=\lambda _n.$ It follows that 
$N_n<N_{n+1}<N$ for all $ n\in \mathbb{N} $ and $p^\bot N_n (\xi )\rightarrow \xi .$ 
Similarly we can find a sequence $(M_n)_{ n\in \mathbb{N} }$ such that $N<M_{n+1}<M_n$ 
for all $ n\in \mathbb{N}$ and $p^\bot (I-M_n)(\eta )\rightarrow \eta $. For every 
 $ n\in \mathbb{N} $ we choose $\omega _n\in H_2$ such that 
$ \| \theta (M_n)\ominus \theta (N_n)(\omega _n) \| =1.$ The operator 
$$s_n=p^\bot N_n(\xi )\otimes 
(\theta (M_n)\ominus \theta (N_n)(\omega _n)) ^*$$ satisfies $s_n=N_ns_n\theta (N_n)^\bot $ and  
so $s_n\in Ball(Y_0).$ Similarly the operator $$t_n=(\theta (M_n)\ominus \theta 
(N_n)(\omega _n)) 
\otimes p^\bot (I-M_n)(\eta )^*$$ satisfies $t_n=\theta (M_n)t_nM_n^\bot $ and so 
$t_n\in Ball(X_0).$ Now we have 
$$s_nt_n=p^\bot N_n(\xi )\otimes p^\bot 
(I-M_n)(\eta )^* $$ which clearly converges in norm to the operator $\xi\otimes  \eta^* .\qquad  \Box$ 

\begin{lemma}\label{2.3} Suppose that $p^\bot \neq 0, \xi \in H_1$ such that 
$\|p^\bot (\xi )\|=1$ and $q$ is the projection onto the space $\overline{p^\bot 
\cl N_1^{\prime 
\prime}\xi }.$ There exists a sequence of finite rank contractions $(r_n)_{n\in \mathbb{N}}\subset 
A$ converging in SOT topology to the projection $q$ such that $r_n=\|\cdot\|-
\lim_{i\in \mathbb{N}}s_i^n t_i^n$ where $s_i^n\in Ball(R_\infty ^{fin}(Y_0)), t_i^n\in 
Ball(C_\infty ^{fin}(X_0))$ for all $i,n \in \mathbb{N}.$ 
\end{lemma}
\textbf{Proof} We define the continuous order preserving map 
$$\phi : \cl N_1p^\bot \rightarrow [0,1], \phi (Np^\bot )=
\|Np^\bot (\xi )\|^2.$$ Since the nest $\cl N_1p^\bot $ is continuous 
$\phi $ is onto $[0,1].$ 
Choose $N_{k,n}p^\bot $ the least element in $\cl N_1p^\bot $ such that $\phi (N_{k,n}p^\bot )=
\frac{k}{2^n}, k=0,1,...,2^n.$ 

We denote $$E_{k,n}=(N_{k,n}\ominus N_{k-1,n})p^\bot ,\;\; \xi _{k,n}=2^{\frac{n}{2}}
E_{k,n}(\xi ), \;\;
r_n=\sum_{k=2}^{2^n}f_{k,n}$$ where $f_{k,n}=\xi_{k-1,n}\otimes  \xi_{k,n}^*. $  

As in  \cite[Lemma 3.9]{dav} we can prove that $\|r_n\|\leq 1$ and the sequence 
$(r_n)_{n\in \mathbb{N}}$ converges in SOT topology to the operator $q.$ 

By the above lemma there exist sequences of rank 1 operators $(s_i^{k,n})_{i\in \mathbb{N}}\subset 
Ball(Y_0), $$(t_i^{k,n})_{i\in \mathbb{N}}\subset Ball(X_0),$ such that 
$s_i^{k,n}t_i^{k,n}\stackrel{\|\cdot\|}{\rightarrow} f_{k,n}, i\rightarrow \infty $ 
for all $k,n.$ We denote $$ s_i^n=(s_i^{2,n}, s_i^{3,n},...,s_i^{2^n,n}),\;\; 
t_i^n=(t_i^{2,n}, t_i^{3,n},...,t_i^{2^n,n})^T$$ and we have 
$r_n=\|\cdot\|-\lim_is_i^nt_i^n.$ Also $$\|s_i^n\|^2=\|\sum_{k=2}^{2^n} s_i^{k,n} (s_i^{k,n})^* 
\|.$$ We may assume that $s_i^{k,n} =E_{k-1,n} s_i^{k,n} $ so 
$$\|s_i^n\|^2=\|\sum_{k=2}^{2^n} E_{k-1,n}  s_i^{k,n} (s_i^{k,n})^*E_{k-1,n} \|.$$ 
Since $\|s_i^{k,n}\|\leq 1$ and the projections  $(E_{k-1,n})_k$ are pairwise orthogonal 
we have $\|s_i^{n}\|\leq 1.$ Similarly we can prove $\|t_i^{n}\|\leq 1.  \qquad  \Box$

\begin{lemma}\label{2.4}  Suppose that $p^\bot \neq 0.$ There exists a net $(g_\lambda )$ 
of finite rank contractions in $A$ converging in SOT topology to $p^\bot $ such that 
$g_\lambda =\|\cdot\|-\lim_{i\in \mathbb{N}}s_i^\lambda t_i^\lambda $ for all $\lambda $ where 
$s_i^\lambda \in Ball(R_\infty ^{fin}(Y_0)), t_i^\lambda \in Ball(C_\infty ^{fin}(X_0))$ 
for all $i\in \mathbb{N}.$ 
\end{lemma}
\textbf{Proof} Using Zorn's Lemma we find a family of vectors $\xi _k: k\in L$ such that the 
projections $q_k$ onto $\overline{p^\bot \cl N_1^{\prime \prime}\xi_k },\;\; k\in L$ are 
pairwise orthogonal and 
they span $p^\bot.$ We assume that $\|p^\bot (\xi _k)\|=1$ for all $k\in L.$ From 
Lemma \ref{2.3} there exist finite rank contractions $(r_n^k)_{n\in \mathbb{N}}$ such that 
$q_k=SOT-\lim_ {n\in \mathbb{N}} r_n^k $ and \linebreak$r_n^k =\|\cdot\|-\lim_{i\in \mathbb{N}}s_i
^{n,k}t_i^{n,k}$ for sequences $$ (s_i^{n,k})_{i\in \mathbb{N}}\subset 
Ball(R_\infty ^{fin}(Y_0)),\;\; (t_i^{n,k})_{i\in \mathbb{N}}\subset Ball(C_\infty ^{fin}(X_0)).$$
 
We define $\cl F=\{F: F\;\;\text{finite\;\;subset\;\;of}\;\; L\}.$ If $F\in \cl F$ 
and $n\in \mathbb{N}$ we define the finite rank contraction $g_{n,F}=\sum_{k\in F}r_n^k.$ 
The family $(g_{n,F})_{n,F}$ indexed by $\mathbb{N}\times \cl F$ is a net. Fix $\xi \in H_1.$

Observe that for all $n\in \mathbb{N}$
$$\| r_n^k (\xi )- q_k(\xi )\|^2=\|q_k(r_n^k -I_{H_1})q_k (\xi )\|^2\leq 
2\|q_k(\xi )\|^2$$
and so $$\sum_{k\in L}\| r_n^k (\xi )-q_k(\xi )\|^2\leq 
2\sum_{k\in L}\|q_k(\xi )\|^2<\infty .$$
If $n\in \mathbb{N}$ and $F\in \cl F$ we have

\begin{equation}\label{ex}  \|g_{n,F}(\xi )-p^\bot (\xi )\|^2 =\|g_{n,F}(\xi )-\sum_{k\in L}
q_k(\xi )\|^2\end{equation} \begin{align*} = &\sum_{k\in F}\|r_n^k(\xi )-
q_k(\xi )\|^2+\|p^\bot (\xi )\|^2-
\sum_{k\in F}\| q_k(\xi )\|^2\\ 
 \leq & \sum_{k\in L} \|r_n^k(\xi )- q_k(\xi )\|^2 +\|p^\bot (\xi )\|^2-
\sum_{k\in F}\| q_k(\xi )\|^2\end{align*}

Since $\lim_{n\in \mathbb{N}}\|r_n^k(\xi )- q_k(\xi )\|^2 =0$ by the Theorem of 
dominated convergence we have $$\lim_{n\in \mathbb{N}}\sum_{k\in L}
\|r_n^k(\xi )- q_k(\xi )\|^2 =0.$$ It follows now from (\ref{ex})
that $\lim_{(n,F)} \|g_{n,F}(\xi )-p^\bot (\xi )\|^2 =0.$ We proved that 
 $SOT-\lim_{(n,F)} g_{n,F}=p^\bot .$ 

If $F=\{k_1,...,k_r\}\subset L$ then $g_{n,F}=
\sum_{m=1}^r r_n^{k_m} $ where $$r_n^{k_m}= \|\cdot\|-\lim_{i\in \mathbb{N}} 
s_i^{n,k_m} t_i^{n,k_m} .$$ So $g_{n,F}=\|\cdot\|-\lim_{i\in \mathbb{N}} s_i^{n,F} t_i^{n,F}$
 where $$s_i^{n,F} =(s_i^{n,k_1},...,s_i^{n,k_r}), \;\; 
t_i^{n,F} =(t_i^{n,k_1},...,t_i^{n,k_r})^T.$$ Since $ s_i^{n,k_j}= q_{k_j} s_i^{n,k_j},\;\; 
t_i^{n,k_j}=t_i^{n,k_j}q_{k_j} $ and the projections $(q_{k_j})$ are pairwise orthogonal 
we conclude that $ s_i^{n,F} \in Ball(R_\infty ^{fin}(Y_0)),  t_i^{n,F} \in Ball(C
_\infty ^{fin}(X_0))$ for all $(n,F).$ This completes the proof. $\qquad \Box$

\begin{theorem}\label{2.6} There exists a net of finite rank  contractions $(f_\lambda )
_{\lambda \in \Lambda }$ converging in SOT topology to the identity operator $I_{H_1}$ 
such that $f_\lambda =\|\cdot\|-\lim_{i\in \mathbb{N}}v_i^\lambda u_i^\lambda $ where 
$(v_i^\lambda)_{i\in \mathbb{N}} \subset Ball(R_\infty ^{fin}(Y_0)), (u_i^\lambda)_{i\in 
\mathbb{N}}\subset Ball(C_\infty ^{fin}(X_0))$ for all $\lambda \in \Lambda.$ 
\end{theorem}
\textbf{Proof} If $p^\bot =0$ the conclusion comes from Lemma \ref{2.1}. So we may assume 
that $p^\bot \neq 0.$ From Lemmas \ref{2.1}, \ref{2.4} there exists a  net 
$(l_\lambda )_{\lambda \in \Lambda} $ of finite rank contractions converging 
 in SOT topology to the projection $p$ such that $l_\lambda =s_\lambda t_\lambda $ 
where $ s_\lambda \in Ball(R_\infty ^{fin}(Y_0)), t_\lambda \in Ball(C_\infty 
^{fin}(X_0))$ for all $\lambda \in \Lambda ,$ and 
a net $(g_\lambda )_{\lambda \in \Lambda }$ 
of finite rank contractions converging in SOT topology to $p^\bot $ such that 
$g_\lambda =\|\cdot\|-\lim_{i\in \mathbb{N}}y_i^\lambda x_i^\lambda $ for all $\lambda \in \Lambda 
$ where $$y_i^\lambda \in Ball(R_\infty ^{fin}(Y_0)),\;\; 
x_i^\lambda \in Ball(C_\infty ^{fin}(X_0))$$ 
for all $i\in \mathbb{N}.$

We denote $f_\lambda =l_\lambda +g_\lambda ,\;\; v_i^\lambda =(s_\lambda ,y_i^\lambda ),\;\;
 u_i^\lambda =(t_\lambda ,x_i^\lambda )^T$ for all $\lambda \in \Lambda , i\in \mathbb{N}.$ 
Observe that $I_{H_1}=SOT-\lim_{\lambda\in \Lambda } f_\lambda $ and 
$f_\lambda =\|\cdot\|-\lim_{i\in \mathbb{N}}v_i^\lambda u_i^\lambda .$ Now we have 
\begin{align*}& \|v_i^\lambda \|^2= \|s_\lambda s_\lambda ^*+y_i^\lambda (y_i^\lambda)^* \| \\
=& \|ps_\lambda s_\lambda ^*p+p^\bot y_i^\lambda (y_i^\lambda)^*p^\bot  \|\leq 1. 
\end{align*} 
Similarly $\|u_i^\lambda\|\leq 1 $ for all $\lambda \in \Lambda , i\in \mathbb{N}.\qquad \Box$

\begin{theorem}\label{2.7} The algebras $A_0, B_0$ are strongly Morita equivalent. 
Particularly $A_0\cong Y_0\otimes^h_{B_0} X_0,$ $B_0\cong X_0\otimes^h_{A_0} Y_0$ 
as operator modules.
\end{theorem} 
\textbf{Proof} We define the bilinear map $Y_0\times X_0\rightarrow A_0: (y,x)\rightarrow yx.$ 
This map is completely contractive and $B_0-$balanced, so induces a completely 
contractive $A_0-$module map $\pi : Y_0\otimes ^h_{B_0} X_0\rightarrow A_0: y \otimes _{B_0} x\rightarrow yx.$ 
We shall prove that $\pi $ is completely isometric.
 It suffices to prove that if $$z_{i,j}= \sum_{k=1}^{m_{i,j}}y_k^{i,j}\otimes _{B_0} x_k^{i,j}, i,j=
1,...,n $$  
then $$\nor{(z_{i,j}) _{i,j} }\leq \nor{\left(\sum_{k=1}^{m_{i,j}}y_k^{i,j} x_k^{i,j}\right)_{i,j} 
 }.$$ 
We recall the contractions $f_\lambda , (v_s^\lambda )_{s\in \mathbb{N}}, 
(u_s^\lambda )_{s\in \mathbb{N}}, \lambda \in \Lambda $ from Theorem \ref{2.6}.

If $x$ is a compact operator then $x=\|\cdot\|-\lim_\lambda xf_\lambda $ (\cite[Proposition
 1.18]{dav}).It follows that 
$ z_{i,j}=\|\cdot\|-\lim_\lambda  \sum_{k=1}^{m_{i,j}}y_k^{i,j}\otimes _{B_0} 
(x_k^{i,j}f_\lambda ). $ If $ \epsilon >0$ there exists $\lambda \in \Lambda $ such that 
$$ \|(z_{i,j})_{i,j} \|-\epsilon  <  \nor{ \left( \sum_{k=1}^{m_{i,j}}y_k^{i,j}\otimes _{B_0} 
(x_k^{i,j}f_\lambda )\right)_{i,j} }-\frac{\epsilon }{2} . $$
Since $ x_k^{i,j} f_\lambda =\|\cdot\|-\lim_{s\in \mathbb{N}}x_k^{i,j}v_s^\lambda u_s^\lambda $
 there exists $s\in \mathbb{N}$ such that 
$$ \nor{  \left(\sum_{k=1}^{m_{i,j}}y_k^{i,j}\otimes _{B_0} 
(x_k^{i,j}f_\lambda )\right)_{i,j} }-\frac{\epsilon }{2} < 
\nor{ \left( \sum_{k=1}^{m_{i,j}}y_k^{i,j}\otimes
 _{B_0} 
( x_k^{i,j}v_s^\lambda  u_s^\lambda )\right)_{i,j} }.$$

Since $x_k^{i,j}v_s^\lambda  \in R_\infty ^{fin}(B_0)$ we have 

\begin{align*}&\|(z_{i,j})_{i,j} \|-\epsilon  <  \nor{\left(\sum_{k=1}^{m_{i,j}}(y_k^{i,j}
x_k^{i,j}
v_s^\lambda)  
\otimes _{B_0}  u_s^\lambda \right)_{i,j} }  \\
= & \nor{\left(\sum_{k=1}^{m_{i,j}}y_k^{i,j}x_k^{i,j}\right)_{i,j} \left( v_s^\lambda  
\otimes _{B_0}  u_s^\lambda \oplus...\oplus  v_s^\lambda  
\otimes _{B_0}  u_s^\lambda \right)} \\ 
\leq  & \nor{\left(\sum_{k=1}^{m_{i,j}}y_k^{i,j}x_k^{i,j}\right)_{i,j} }\nor{ v_s^\lambda  
\otimes _{B_0}  u_s^\lambda}
\\ \leq & 
\nor{\left(\sum_{k=1}^{m_{i,j}}y_k^{i,j}x_k^{i,j}\right)_{i,j} } \nor{ v_s^\lambda}\nor{  
u_s^\lambda} \leq \nor{\left(\sum_{k=1}^{m_{i,j}}y_k^{i,j}x_k^{i,j}\right)_{i,j} } 
\end{align*}

Since $\epsilon $ was arbitrary we have 
$\nor{(z_{i,j})_{i,j} }\leq \nor{\left(\sum_{k=1}^{m_{i,j}}y_k^{i,j}\otimes _{B_0} x_k^{i,j}
\right)_{i,j}  }.$
We proved that $\pi $ is completely isometric. It remains to prove that $\pi $ is 
onto $A_0.$ It suffices to prove that the space $Im\pi $ is dense in $A_0.$ 

Let $a\in Ball(A_0)$ and $\epsilon >0.$ Since $a= \|\cdot\|-\lim_ \lambda f_\lambda a$ there 
exists $\lambda \in \Lambda $ such that $\|a-f_\lambda a\| <\frac{\epsilon }{2} .$ 
Since $f_\lambda =\|\cdot\|-\lim_s v_s^\lambda u_s^\lambda  $ there 
exists $s\in \mathbb{N}$ such that 
$$ \|f_\lambda -v_s^\lambda u_s^\lambda  \|<\frac{\epsilon }{2} .$$
It follows that $\|a - v_s^\lambda u_s^\lambda a  \|< \epsilon .$ 
But $v_s^\lambda u_s^\lambda a=\pi (v_s^\lambda \otimes _{B_0}(u_s^\lambda a))$ 
and this completes the proof. Similarly we can prove that $B_0\cong X_0\otimes^h_{A_0} Y_0.
\qquad \Box$

\bigskip

We define the bilinear map $Y\times X\rightarrow A: (y,x)\rightarrow yx.$ 
This map is completely contractive $B-$balanced and separately $w^*-$continuous, 
so induces a completely 
contractive $w^*-$continuous map $\rho : Y\otimes ^{\sigma h}_{B} X\rightarrow A: 
y \otimes _{B} x\rightarrow yx$ which is also an $A-$module map.  
We shall prove that the restriction of $\rho $ on the space 
$Y\otimes ^{h}_{B} X$ is completely isometric and we shall use this fact in Theorem \ref{2.11} to
 prove that $A\cong Y\otimes ^{\sigma h}_BX.$ 

\begin{lemma}\label{2.8} The restriction of $\rho $ on the space 
$Y\otimes ^{h}_{B} X$ is completely isometric.
\end{lemma}
\textbf{Proof} It suffices to prove that if 
$$z_{i,j}= \sum_{k=1}^{m_{i,j}}y_k^{i,j}\otimes _{B} x_k^{i,j}, i,j=
1,...,n $$  
then $$\nor{(z_{i,j})_{i,j} }\leq \nor{\left(\sum_{k=1}^{m_{i,j}}y_k^{i,j} x_k^{i,j}\right)_{i,j}  }.$$ 
We recall the contractions $f_\lambda , (v_s^\lambda )_{s\in \mathbb{N}}, 
(u_s^\lambda )_{s\in \mathbb{N}}, \lambda \in \Lambda $ from Theorem \ref{2.6}.
Fix $\lambda \in \Lambda .$ If $\epsilon >0$ there exists $s\in \mathbb{N}$ such that 
\begin{align*} & \nor{\left(\sum_{k=1}^{m_{i,j}} y_k^{i,j} \otimes _B(x_k^{i,j} f_\lambda )
\right
)_{i,j} }
-\epsilon < \nor{\left(\sum_{k=1}^{m_{i,j}} y_k^{i,j} \otimes _B(x_k^{i,j} 
v_s^\lambda u_s^\lambda 
)\right)_{i,j} } \\ =& \nor{\left(\sum_{k=1}^{m_{i,j}} (y_k^{i,j}x_k^{i,j} v_s^\lambda )\otimes 
_B u_s^\lambda 
\right)_{i,j} }=  \nor{\left(\sum_{k=1}^{m_{i,j}} y_k^{i,j}x_k^{i,j}\right)_{i,j} (  v_s^\lambda
 \otimes _B u_s^\lambda 
\oplus ...\oplus v_s^\lambda \otimes _B u_s^\lambda )}\\\leq & 
\nor{\left(\sum_{k=1}^{m_{i,j}}y_k^{i,j} x_k^{i,j}\right)_{i,j}  }
\end{align*}
It follows that $$\nor{ \left(\sum_{k=1}^{m_{i,j}} y_k^{i,j} \otimes _B(x_k^{i,j} f_\lambda )
\right)_{i,j}  }\leq 
\nor{\left(\sum_{k=1}^{m_{i,j}}y_k^{i,j} x_k^{i,j}\right)_{i,j}  }$$ for all $\lambda \in \Lambda .$
Since  $$\left(\sum_{k=1}^{m_{i,j}} y_k^{i,j} \otimes _B x_k^{i,j} \right)_{i,j} = 
w^*-\lim_\lambda \left(\sum_{k=1}^{m_{i,j}} y_k^{i,j} \otimes _B(x_k^{i,j} f_\lambda )\right)_{i,j}  $$
 we have  $$\nor{\left(\sum_{k=1}^{m_{i,j}} y_k^{i,j} \otimes _B x_k^{i,j} \right)_{i,j} }\leq  
\nor{\left(\sum_{k=1}^{m_{i,j}} y_k^{i,j} x_k^{i,j} \right)_{i,j} }.\qquad \Box $$

\bigskip

The second dual operator space $A_0^{**}$ of the operator algebra $A_0$ is also an 
operator algebra with product describing in \cite[section 2.5]{bm}. The product on 
$A_0^{**}$ extends the product on $A_0.$ With this we mean that if $\iota : A_0\rightarrow 
A_0^{**}$ is the canonical embedding then $\iota(ab)= \iota(a) \iota(b) $ 
for all $a,b\in A_0.$

\begin{lemma}\label{2.9} The operator algebra $A, $ (resp. $B$) is isomorphic  
as dual operator algebra with $A_0^{**}$ (resp. $B_0^{**}$ ).
\end{lemma}
\textbf{Proof} We denote by $C_1$ the space of trace class operators in $H_1$ and 
$$ \Omega =\{c\in C_1: (N_-)^\bot cN=0\;\; \forall N\;\;\in\;\; \cl N_1\}.$$
By \cite[section 16]{dav} the maps 
$$\mu : C_1/ \Omega \rightarrow A_0^*: \mu (c)(a)=tr(ca),$$
$$\sigma : A\rightarrow (C_1/ \Omega)^* : \sigma (a)(c)=tr(ac)$$
are surjective isometries. We define the isometry $\phi =(\mu ^*)^{-1}\circ \sigma :
A\rightarrow A_0^{**}.$ This map satisfies $\phi (a)=\iota (a)$ for all $a\in A_0.$ 
Since $i(ab)=\iota(a) \iota (b)$ for all $a,b \in A_0$ and $\phi $ is $w^*-$continuous 
$\phi $ is a homomorphism onto $A_0^{**}.$ (When we say $w^*-$continuous we mean 
that $\phi $ is $B(H_1)_*-A_0^*$ continuous.)

If $n\in \mathbb{N}$ the algebra $M_n(A)$ is also a nest algebra, so by the above arguments, 
there exists a  $w^*-$continuous isometry $$ \stackrel{\sim  }{\phi }: M_n(A)\rightarrow 
(M_n(A)_0)^{**}= M_n(A_0)^{**} $$ such that $$\stackrel{\sim  }{\phi }((a_{i,j}))=
 \stackrel{\sim  }{\iota  }((a_{i,j}))$$ for all $(a_{i,j})\in M_n(A_0),$ 
where $\stackrel{\sim  }{\iota }: M_n(A_0)\rightarrow M_n(A_0)^{**}  $ is the 
canonical embedding. By \cite[1.4.11]{bm} there exists a $w^*-$continuous 
isometry $\tau : M_n(A_0)^{**}  \rightarrow M_n(A_0^{**})$ such that 
$\tau (\stackrel{\sim  }{\iota  }((a_{i,j}))=(\iota (a_{i,j}))$ for all 
$(a_{i,j})\in M_n(A_0).$ So we have a $w^*-$continuous 
isometry $ \tau \circ \stackrel{\sim  }{\phi } : M_n(A)\rightarrow M_n(A_0^{**}) $ 
satisfying $$\tau \circ \stackrel{\sim  }{\phi } ((a_{i,j}))=(\iota (a_{i,j}))$$ 
for all $(a_{i,j})\in M_n(A_0).$ But the map 
$$\phi _n: M_n(A)\rightarrow M_n(A_0^{**}) : (b_{i,j})= (\phi (b_{i,j})) $$  is 
a  $w^*-$continuous map satisfying $$\phi _n((a_{i,j}))=(\phi (a_{i,j}))=(\iota (a_{i,j}))  $$
 for all $(a_{i,j})\in M_n(A_0).$ So $\phi _n$ is equal to $ \tau \circ \stackrel{\sim  }{\phi } 
$ in $M_n(A_0).$ Since $\overline{M_n(A_0)}^{w^*}=M_n(A)$ we have 
$\phi _n=\tau \circ \stackrel{\sim  }{\phi }. $ So $\phi _n$ is isometry for all $n\in \mathbb{N}
.$ We proved that $\phi $ is a completely isometric map and this completes the proof.$\qquad \Box$

\bigskip

We are now ready to present the main theorem of this paper:

\begin{theorem}\label{2.11} A. The following are equivalent:

(i) The nests $\cl N_1, \cl N_2$ are isomorphic.

(ii) The algebras $A_0, B_0$ are strongly Morita equivalent.

(iii) The algebras $A, B$ are weakly$-*$ Morita equivalent.

B. If $\theta : \cl N_1\rightarrow \cl N_2$ is a nest isomorphism, $X=Op(\theta ), Y=Op(\theta 
^{-1})$ then:

(i)  $ A_0\cong Y_0\otimes ^h_{B_0}X_0,\;\;  B_0\cong X_0\otimes ^h_{A_0}Y_0,$ as operator modules,

(ii) $A\cong Y\otimes ^{\sigma h}_BX,\;\;  B\cong X\otimes ^{\sigma h}_AY,$ as dual operator modules.
\end{theorem}

\textbf{Proof}

A.(i)$\Rightarrow$ (ii)

This is Theorem \ref{2.7}.

(ii)$\Rightarrow $(iii)

 If $A_0$ and $B_0$ are strongly Morita equivalent 
then the operator algebras 
$A_0^{**}$ and $B_0^{**}$ are weakly$-*$ Morita equivalent, \cite[section 3]{bk} . So by Lemma \ref{2.9} $A$ and $B$ 
are weakly$-*$ Morita equivalent.

(iii)$\Rightarrow $(iv)

Let $(A, B, V, U)$ be a weak$-*$ Morita context \cite{bk}. It follows that there exist 
completely contractive separately $w^*-$continuous bilinear maps $(\cdot,\cdot): 
V\times U \rightarrow A$ which is $A-$module and $B-$balanced map and  $[\cdot,\cdot]: 
U\times V \rightarrow B$ which is $B-$module and $A-$balanced map satisfying 
$$(y,x)y^\prime=y[x,y^\prime],\;\;  x^\prime(y,x)=[x^\prime,y]x\;\;  \forall x,x^\prime\in 
U,\;\; y,y^\prime\in V,$$
and $ A=\wsp{\{(y,x): x\in U, y\in V\}},  B=\wsp{\{[x,y]: x\in U, y\in V\}}.$ 

If $N\in \cl N_1$ we define $\theta (N)$ the projection onto the space generated  
by vectors of the form $[xN,y](\omega ), x\in U, y\in V, \omega \in H_2.$ Since 
$b[xN,y]=[bxN,y]$ for all $b\in B$ we have $\theta (N)^\bot B\theta (N)=0$ so $\theta (N)\in \cl N_2.$ 
Also if $N_1\leq N_2$ then $\theta(N_1)\leq  \theta (N_2)$ and so $\theta $ is 
an order preserving map from $\cl N_1$ into $\cl N_2.$ 

Similarly if $M\in \cl N_1$ we define $\sigma (M)$ the projection onto the space generated  
by vectors of the form $(yM,x)(\xi ), x\in U, y\in V, \xi \in H_1.$ The 
map $\sigma : \cl N_2\rightarrow \cl N_1$  is 
an order preserving map.

If $x, x^\prime\in U, y\in V$ and $N\in \cl N_1$ then 
\begin{align*}& \theta (N)^\bot [xN,y]=0 \Rightarrow [\theta (N)^\bot xN,y]=0 \Rightarrow 
[\theta (N)^\bot xN,y]x^\prime=0\\ \Rightarrow & \theta (N)^\bot xN (y,x^\prime)=0. \end{align*}
Since the operators $(y,x^\prime)$ span the algebra $A$ we have 
\begin{equation}\label{ex1}\theta (N)^\bot xN =0\Rightarrow xN=\theta (N)xN
\;\; \forall x\in U,\;\; N\in 
\cl N_1.
\end{equation}

Similarly 
\begin{equation}\label{ex2} yM=\sigma (M)yM \;\;\forall\;\; y\in V,\;\; M\in \cl N_2.
\end{equation}
   
If $x, x^\prime \in U, y, y^\prime\in V, N\in \cl N_1$ we have 
\begin{align*}& [xN^\bot ,y] [x^\prime N, y^\prime ]= [[xN^\bot ,y] x^\prime N,y^\prime ] \\
=&[xN^\bot (y, x^\prime) N,y^\prime ]= 0 \;\;\text{because}\;\;(y, x^\prime) \in 
\mathrm{Alg}(\cl N_1)
\end{align*}

It follows that $[x N^\bot ,y]\theta (N)=0\Rightarrow [x, N^\bot y\theta (N) ]= 0$ 
for all $x\in U, y\in V, $ and so 

\begin{equation}\label{ex3} N^\bot y\theta (N) =0\Rightarrow y\theta (N)=Ny\theta (N),\;\;
 \forall\;\; 
y\in V,\;\; N\in \cl N_1 
\end{equation}
 
 Similarly we can prove 

\begin{equation}\label{ex4} x\sigma (M)=Mx\sigma (M)\;\; \forall \;\;x\in U, \;\;M\in \cl N_2
\end{equation}

If $N\in \cl N_1$ and $x\in U, y\in V$ then 
$(y,x)N=(y,xN)=(y,\theta (N)xN)$ because of (\ref{ex1}). The last operator is equal to 
$(y\theta (N),xN)=\sigma (\theta (N))(y,x)N$ because of (\ref{ex2}).
It follows that $N\leq \sigma (\theta (N)).$

Similarly $(y,x)^*N^\bot =(N^\bot (y,x))^*=(N^\bot y,x)^*= (N^\bot y\theta (N)^\bot ,x)^* $ 
because of (\ref{ex3}). The last operator is equal to $$ (N^\bot y, \theta (N)^\bot x)^* =
(N^\bot y, \theta (N)^\bot x\sigma( \theta(N))^\bot  )^* $$ because of (\ref{ex4}). 
The last operator is equal to $\sigma( \theta (N))^\bot (N^\bot y,\theta (N)^\bot x)^*.$

Since $I_{H_1}=w^*-\lim_i\sum_{k=1}^{n_i}( y_k^i , x_k^i )^*$ for $(y_k^i)\subset V, (x_k^i )\subset U
$ we have $N^\bot \leq \sigma (\theta (N))^\bot $ and so $N= \sigma (\theta (N)).$

Similarly we can prove that $M=\theta (\sigma (M))$ for all $M\in \cl N_2.$ This 
completes the proof of the fact that $\theta $ is a nest isomorphism.

\bigskip

B. Let $\theta : \cl N_1\rightarrow \cl N_2$ be a nest isomorphism and  $X=Op(\theta ), 
Y=Op(\theta^{-1}).$  Claim B-(i) follows from Theorem \ref{2.7}.

Let $\rho : Y\otimes ^{\sigma h}_{B} X \rightarrow A $ be the map which was defined above Lemma 
\ref{2.8}. Let $z\in Ball(M_n(Y\otimes ^{\sigma h}_{B} X))$ for a fixed $n\in \mathbb{N}.$
By \cite[Corollary 2.8]{bk} there exists a net $(z_i)\subset Ball(M_n( Y\otimes ^{h}_{B}X ))$ 
 converging in $w^*$ topology to $z.$ It follows that $\rho (z_i)\stackrel{w^*}
{\rightarrow}\rho (z)$ in $M_n(A).$ If $f,g$ are finite rank operators in $A$ we denote 
$f^n=f\oplus f\oplus ...\oplus f$ and similarly for $g^n.$ We have that 
$$ f^n\rho (z_i)g^n\stackrel{\|\cdot\|}{\rightarrow}f^n\rho (z)g^n \Rightarrow  
\rho (f^n z_i g^n) \stackrel{\|\cdot\|}{\rightarrow}\rho ( f^n zg^n ). $$

From Lemma \ref{2.8} it follows that 
$ f^n zg^n \in M_n(Y\otimes ^{h}_{B} X)$ and $\|\rho (f^n zg^n )\|=\|f^n zg^n \|$ 
for all finite rank operators $f, g$ in $A.$ We recall the finite rank contractions 
$(f_\lambda )_{\lambda \in \Lambda} $ from Theorem \ref{2.6}. 
For $\lambda, \mu \in \Lambda $ we have 
$$ \| f^n_\lambda zf^n_\mu \|= \|\rho (f^n_\lambda zf^n_\mu)\|= \|f^n_\lambda \rho (z)
f^n_\mu\|\leq \|\rho (z)\|.$$

Since $ zf^n_\mu =w^*-\lim_\lambda f^n_\lambda zf^n_\mu $ we have $\| zf^n_\mu  \|\leq \|
\rho (z)\|$ for all $\mu \in \Lambda .$ Now taking the $w^*-$limit of $(zf^n_\mu)_{\mu
\in \Lambda}  $
 we obtain $\|z\|\leq \|\rho (z)\|.$ We proved that the map 
$\rho : Y\otimes ^{\sigma h}_{B} X \rightarrow A $ is a complete isometry. From Theorem 
\ref{2.7} and its proof we have that $A_0=\lsp{Y_0X_0}.$ 
Since $A=\overline{A_0}^{w^*}$ we have 
$$A=\wsp{YX}=\wsp{\{\rho (y\otimes _Bx): y\in Y, x\in X\}}.$$
 By the Krein-Smulian Theorem the space $Im\rho $ is $w^*-$closed and so $\rho $ is onto $A.$ 
Similarly we can prove that $B\cong X\otimes ^{\sigma h}_AY,$ as dual operator modules.
$\qquad \Box$

\section{Spatial Morita equivalence and nest algebras}

In this section we shall investigate the relation between weak$-*$ and 
spatial Morita equivalence for nest algebras. We give the definition of spatial 
Morita equivalence:

\begin{definition}\label{spat1} (I. G. Todorov) Let $C, D$ be $w^*-$closed algebras acting on 
the Hilbert spaces $K_1, K_2$ respectively. 
 We say that $C$ and $D$ are \textbf{spatially Morita equivalent} if there exists a
 $D-C$ bimodule $V\subset B(K_1,K_2)$ and  a $C-D$ bimodule $U\subset B(K_2, K_1)$ 
such that $C=\wsp{UV}, D=\wsp{VU}.$ 
\end{definition}

We also need the following notions. 
If $\cl L$ is a set of projections acting on the Hilbert space $H$ 
the set $$ \mathrm{Alg}(\cl L) = \{x\in B(H): p^\bot xp=0,\;\; \forall\;\; p\in \cl L\} $$ 
is an algebra. An algebra $A$ is called reflexive if there exists a  set of projections 
$\cl L$ such that $A=\mathrm{Alg}(\cl L).$ In the special case where $\cl L$ is a 
complete lattice of commuting projections containing the zero and identity operators 
 the algebra $\mathrm{Alg}(\cl L)$ is 
called a CSL algebra and the lattice $\cl L$ is called a CSL lattice. 
Obviously, nest algebras are CSL algebras. If $A$ is an 
algebra acting on the Hilbert space $H$ the lattice  
$$\{p\in pr(B(H)): p^\bot xp=0,\;\; \forall\;\; x\in A\}$$ is called the 
lattice of $A$ and we denote it by $\mathrm{Lat}(A).$ If $\cl L$ 
is a CSL lattice then $\mathrm{Lat}(\mathrm{Alg}(\cl L))=\cl L,$ \cite{arv}, \cite{dav2}.  

Two spatially Morita equivalent algebras 
are not always weakly$-*$ Morita equivalent even in the case one of them is a CSL 
algebra:

\begin{example}\label{spat2}\em{Let $C$ be a nest algebra. We denote the algebras 
$A=C\oplus C$ and $$ B=\{\left(\begin{array}{clr} a & b-a \\ 0 & b \end{array}\right): 
a, b \in C\}
.$$ Observe that $A$ is a CSL algebra whose lattice is 
$$ \mathrm{Lat}(A) =\{p\oplus q: p,q \in  \mathrm{Lat}(C) \}.$$ 
Since the center of $C$ is trivial \cite[Corollary 19.5]{dav} the center of 
$A$ is $Z(A)= \mathbb{C} \oplus \mathbb{C} $ and the center 
of $B$ is $$Z(B)=
\{\left(\begin{array}{clr} \lambda  & \mu- \lambda  \\ 0 & \mu  \end{array}\right): \lambda, 
\mu  \in \mathbb{C}\}.$$ 
We also denote the spaces 
$$X=\{\left(\begin{array}{clr} a & -a \\ 0 & \;\;\;b \end{array}\right): 
a, b \in C\},\;\;\; Y= \{\left(\begin{array}{clr} a & b \\ 0 & b \end{array}\right) : 
a, b \in C\}.$$
We can check that $X$ is an $A-B$ bimodule, $Y$ is a $B-A$ bimodule and $XY=A,\;\;\; YX=B.$ 
So the algebras $A, B$ are spatially Morita equivalent. If $A$ and $B$ were 
weakly$-*$ Morita equivalent by \cite[Theorem 3.7]{bk} they would have isomorphic centers 
through a completely isometric homomorphism. This is a contradiction because 
$Z(A)$ is a von Neumann algebra and $Z(B)$ is a non-selfadjoint algebra.}
\end{example}

Despite the above example, in \cite{ele} we proved that two CSL algebras 
are spatially Morita equivalent if and only if their lattices are isomorphic, so by Theorem \ref{2.11} 
we conclude the following theorem:

\begin{theorem}\label{spat3} Two nest algebras are spatially Morita equivalent if 
and only if they are 
weakly$-*$ Morita equivalent.
\end{theorem}

Also despite the example \ref{spat2} we have the following theorem:

\begin{theorem}\label{spat4} Let $A$ be a nest algebra, $B$ be a unital dual operator algebra 
and $\beta $ be a completely isometric normal representation of $B$ such that 
$A$ and $\beta (B)$ are spatially Morita equivalent. It follows that $A$ and $B$ are 
weakly$-*$ Morita equivalent.
\end{theorem}
\textbf{Proof} By \cite[Theorem 4.1, remark 4.2]{ele} $\beta (B)$ is a nest algebra 
whose nest is isomorphic with the nest of $A.$ The conclusion comes 
from Theorems \ref{2.11} and \ref{spat3}.$\qquad \Box$

\begin{theorem}\label{spat5}(Blecher-Kashyap) If $A, B$ are weakly$-*$ Morita equivalent
 unital dual operator algebras, for every completely isometric normal representation 
$\alpha $ of $A$ there exists a completely isometric normal representation 
$\beta $ of $B$ such that the algebras $\alpha (A), \beta (B)$ are 
spatially Morita equivalent.   
\end{theorem}
\textbf{Proof} Suppose that $(A, B, X, Y)$ is a weakly$-*$ Morita context \cite{bk}. We use now 
 arguments from the beginning of the 4th section of \cite{bk}. If 
$\alpha $ is a completely isometric normal representation of $A$ on the Hilbert space 
$H$ the tensor product $K=Y\otimes ^{\sigma h}_A H$ with its norm 
is a Hilbert space on which $B$ is represented through the $w^*-$continuous  
complete isometry $\beta $ given by 
$$\beta (b)(y\otimes h)=(by)\otimes h\;\; \forall \;\;b\in B,\;\; y\in Y,\;\; h\in H.$$
 Also Blecher and Kashyap prove that the maps $\phi: Y\rightarrow B(H, K), \psi: X\rightarrow 
B(K, H) $ given by 
$\phi (y)(h)=y\otimes h$ and $\psi (x)(y\otimes h)=\alpha ((x,y))(h)$ are 
$w^*-$continuous complete isometries. See in \cite{bk} for the properties 
of the bilinear map $(\cdot, \cdot): X\times Y\rightarrow A.$ We can easily check that $\psi (X)$
 is an $\alpha(A)- \beta(B) $ bimodule, $\phi (Y)$ is a $\beta(B) -\alpha(A) $ 
bimodule and $$\alpha (A)=\wsp{\psi (X)\phi (Y)},\;\;\; 
\beta (B)=\wsp{\phi (Y)\psi (X)} .\qquad \Box$$ 

\begin{corollary}\label{spat6} If $A$ is a nest algebra and $B$ is a unital dual operator algebra 
which are weakly$-*$ Morita equivalent then there exists a completely isometric 
normal representation $\beta $ of $B$ such that $\beta (B)$ is a nest algebra.
\end{corollary}
\textbf{Proof} By  the above Theorem  there exists a completely isometric 
normal representation $\beta $ of $B$ such that the algebras $A$ and 
$\beta (B)$ are spatially Morita equivalent. From \cite[remark 4.1]{ele} 
the algebra $\beta (B)$ is reflexive and from \cite[Theorem 4.2]{ele} the lattice of 
$\beta (B)$ is isomorphic with the nest of $A.$ So $\beta (B)$ is a nest algebra. $\qquad \Box$

\begin{corollary}\label{spat7} If $A$ is a CSL algebra which is not a nest algebra 
then $A$ is not weakly$-*$ Morita equivalent with anyone nest algebra.
\end{corollary}  
\textbf{Proof} By the above corollary if $A$ was weakly$-*$ Morita equivalent with a nest algebra then 
it would have a normal completely isometric representation $\alpha $ 
such that $\alpha (A)$ is a nest algebra. This is a contradiction because as we can easily check 
$\alpha (\mathrm{Lat}(A))=\mathrm{Lat}(\alpha (A)). \qquad \Box$

\section{A stable isomorphism theorem for nest algebras}

In this section we are going to present a new theorem which characterizes the stable isomorphism 
of separably acting nest algebras.

\begin{definition}\label{stable} Two dual operator algebras $C, D$ are called 
\textbf{stably isomorphic} if there exists a Hilbert space $H$ and a 
completely isometric, $w^*$-bicontinuous isomorphism from the algebra $C \stackrel{-}{\otimes} B(H)$ 
onto the algebra $D \stackrel{-}{\otimes} B(H),$ 
where $\stackrel{-}{\otimes} $ is the normal spatial tensor product. 
\end{definition}

We give two relevant definitions:

\begin{definition}\label{3.1} \cite{ele} Let $C, D$ be $w^*$ closed algebras 
acting on Hilbert spaces $H_1$ and $H_2$ respectively. If there exists a TRO 
$\cl{M}\subset B(H_1,H_2),$ i.e. a subspace satisfying $\cl{MM^*M}\subset \cl M,$ 
such that $C=\wsp{\cl{M}^*D\cl{M}}\;\; 
\text{and}\;\;
D=\wsp{\cl{M}C\cl{M}^*}$ we write $C 
\stackrel{\cl{M}}{\sim}D.$ We say that the algebras $C, D$
 are \textbf{TRO equivalent} if there exists a TRO 
$\cl{M}$ such that $C 
\stackrel{\cl{M}}{\sim} D.$ 

\end{definition}

\begin{definition}\label{3.2}\cite{ele1} Let $C, D$ be abstract dual operator algebras. 
These algebras are called \textbf{$\Delta -$equivalent} if they have completely isometric normal 
representations $\phi, \psi $ such that the algebras $\phi(C), \psi(D) $ are TRO-equivalent. 
\end{definition}

In \cite{elepaul} we proved the following theorem:

\begin{theorem}\label{3.3} Two unital dual operator algebras are stably isomorphic if and only if 
they are  $\Delta -$equivalent.   
\end{theorem}
  
$\Delta -$equivalence implies weak$-*$ Morita equivalence \cite[section 3]{bk}. The 
converse does not hold. The counterexample is \cite[example 3.7]{ele2}. We shall 
give a new proof of this fact in Theorem \ref{3.8}.

\cite[Theorem 3.2]{ele2} implies the following corollary:

\begin{corollary}\label{3.4} Two nest algebras are $\Delta -$equivalent if and only if 
they are TRO-equivalent.
\end{corollary}   
In what follows if $X$ is a subset of $B(H)$ where $H$ is a Hilbert space 
we denote by $X^\prime$ the commutant of $X$ and by $X^{\prime \prime}$ 
the algebra $(X^\prime)^\prime.$ 
In \cite{ele} we proved the following criterion of TRO-equivalence for reflexive algebras:

\begin{theorem}\label{3.5} Two reflexive algebras $C, D$ are TRO-equivalent 
if and only if there exists a $*-$isomorphism $\delta :(C\cap C^*)^\prime\rightarrow 
(D\cap D^*)^\prime$ such that $\delta (\Lat{C})=\Lat{D}.$ 
\end{theorem}
  
Comparing Theorems \ref{3.3}, \ref{3.5} and Corollary \ref{3.4} we take the following:

\begin{corollary}\label{3.6} The nest algebras $ \Alg{\cl N_1}, \Alg{\cl N_2}$ are 
stably isomorphic if and only if there exists a $*-$isomorphism 
$\delta : \cl N_1^{\prime\prime} \rightarrow \cl N_2^{\prime\prime}$ 
 such that $\delta (\cl N_1)=\cl N_2.$ 

\end{corollary}

In the rest of this section we fix two nests $\cl N_1, \cl N_2$ acting on the 
separable Hilbert spaces $H_1, H_2$ respectively and we denote  $ A=\Alg{\cl N_1}, 
B=\Alg{\cl N_2}.$ We use now extensively notions from \cite[section 7]{dav}. 
If $\xi $ (resp. $\omega $) is a unit separating vector for the algebra 
 $\cl N_1^{\prime\prime} $ (resp. $\cl N_2^{\prime\prime})$ we define the order isomorphism 
$\phi _\xi $ (resp. $\psi _\omega $) from $\cl N_1$ (resp. $\cl N_2$) onto a closed 
 subset of the interval $[0,1]$ given by $\phi _\xi (N)=\|N(\xi )\|^2$ (resp. 
$ \psi _\omega (M)=\|M(\omega )\|^2$). 

Suppose that $[0,1]\setminus \phi_ \xi (\cl N_1) =\cup _{n}(l_n, r_n)$ and 
$[0,1]\setminus \psi_ \omega(\cl N_2) =\cup _n(t_n, s_n). $ If $m$ is the 
Lebesgue measure  we define the  measures $ \mu _\xi , \nu _\omega $ 
given by 
$$\mu_ \xi(S)=m(S\cap \phi_ \xi (\cl N_1) ) +\sum_{r_n\in S}(r_n-l_n)$$
$$\nu _\omega (S)=m(S\cap \psi_ \omega(\cl N_2) )+\sum_{s_n\in S}(s_n-t_n), $$ 
for every Borel subset $S$ of $[0,1].$
We denote $\cl M_1$ (resp. $\cl M_2$) the nest $\{M_s: 0\leq s\leq 1 \}\subset 
B( L^2([0,1], \mu _\xi ) )$ (resp. $\{N_s: 0\leq s\leq 1 \}\subset B(L^2([0,1], \nu _\omega )))
 $ where $M_s$ (resp. $N_s$) is the projection onto the space $ L^2([0,s], \mu _\xi ) $
 (resp. $L^2([0,s], \nu _\omega ) $).

The algebra 
 $\cl N_1^{\prime\prime}$ is $*-$isomorphic with the algebra $ L^\infty ([0,1], \mu_ \xi ) $ 
(resp. $ L^\infty ([0,1], \nu_ \omega ) $) acting on the Hilbert space 
$ L^2([0,1], \mu _\xi ) $ (resp. $ L^2([0,1], \nu _\omega ) $) through an isomorphism 
mapping the nest $\cl N_1$ (resp. $\cl N_2$) onto $\cl M_1$ (resp. $\cl M_2$).

We denote by $AbsHom([0,1])$ the set of order homeomorphisms 
 $\alpha : [0,1]\rightarrow [0,1]$ which satisfy the 
property $m(S)=0\Rightarrow m(\alpha (S))=0.$ The theorem 
below describes when two separably acting nest algebras are stably isomorphic.

\begin{theorem}\label{3.7} The algebras $A, B$ are stably isomorphic if and only if 
there exist separating unit vectors $\xi $ for $\cl N_1^{\prime\prime} ,$  
$\omega $ for $\cl N_2^{\prime\prime} $ and $\alpha \in AbsHom([0,1])$ such that
 $\alpha (\phi_ \xi (\cl N_1) )=\psi_ \omega(\cl N_2) .$ 
\end{theorem}
\textbf{Proof} Suppose that the algebras $A, B$ are stably isomorphic. From Corollary \ref{3.6} 
there exists a $*-$isomorphism 
 $\delta : \cl N_1^{\prime\prime} \rightarrow \cl N_2^{\prime\prime}$ 
 such that $\delta (\cl N_1)=\cl N_2.$ Fix separating unit vectors $\xi $ for 
$\cl N_1^{\prime\prime} ,$ and  
$\omega $ for $\cl N_2^{\prime\prime}. $ Taking compositions we obtain a $*-$isomorphism 
$$ \stackrel{\sim }{\delta } : L^\infty ([0,1], \mu_ \xi ) \rightarrow 
L^\infty ([0,1], \nu_ \omega ) $$ such that $\stackrel{\sim }{\delta } (\cl M_1)=\cl M_2.$ 
Every isomorphism between maximal abelian selfadjoint algebras is implementing by a unitary. 
So the nests $\cl M_1, \cl M_2$ are unitarily equivalent. By \cite[Theorem 7.23]{dav} 
there exists $\alpha \in AbsHom([0,1])$ such that 
$\alpha (\phi_ \xi (\cl N_1) )=\psi_ \omega(\cl N_2) .$ 
 
Conversely if there exist such $\xi, \omega $ and $\alpha ,$ by the same theorem there exists a unitary 
$u\in B(L^2([0,1], \mu _\xi ) ,L^2([0,1], \nu _\omega ) )$ such that $u^*\cl M_2u=\cl M_1.$ 
It follows that  $L^\infty ([0,1], \mu _\xi )= u^*L^\infty ([0,1], \nu _\omega )u. $ 
Taking compositions we take a $*-$isomorphism  
$\delta : \cl N_1^{\prime\prime} \rightarrow \cl N_2^{\prime\prime}$ 
 such that $\delta (\cl N_1)=\cl N_2.$ Again from Corollary \ref{3.6} we conclude that the algebras 
$A$ and $B$ are stably isomorphic. $\qquad \Box$

\begin{remark}\label{separate}\em{If 
there exist separating unit vectors $\xi $ for $\cl N_1^{\prime\prime} ,$  
$\omega $ for $\cl N_2^{\prime\prime} $ and $\alpha \in AbsHom([0,1])$ such that 
$\alpha (\phi_ \xi (\cl N_1) )=\psi_ \omega(\cl N_2) $ then for all 
separating unit vectors 
 $\xi_1 $ for $\cl N_1^{\prime\prime} $ and   
$\omega_1 $ for $\cl N_2^{\prime\prime} $ there exists $\alpha_1 \in AbsHom([0,1])$ such that 
$\alpha _1 (\phi_ {\xi _1} (\cl N_1) )=\psi _ {\omega _1}(\cl N_2) .$ This is a consequence of  
\cite[Proposition 7.22]{dav}.}
\end{remark}

\bigskip

We give a new proof of the following result:

\begin{theorem}\label{3.8} Weak$-*$ Morita equivalence is strictly weaker than $\Delta -$equivalence.
\end{theorem}
\textbf{Proof} Let $C$ be the Cantor set, $\gamma $ be an order homeomorphism of $[0,1]$ 
such that $m(\gamma (C))>0.$ Suppose that  
$[0,1]\setminus C = \cup _{n}(l_n, r_n)$ and 
$[0,1]\setminus   \gamma (C)= \cup _n (t_n, s_n). $ We denote by $\mu $ the measure 
$$\mu(S)=\sum_{r_n\in S}(r_n-l_n)$$ and by $\nu $ the measure 
$$\nu (S)=m(S\cap \gamma (C) )+\sum_{s_n\in S}(s_n-t_n). $$
We denote $\cl M_1$ (resp. $\cl M_2$) the nest $\{M_s: 0\leq s\leq 1 \}\subset 
B( L^2([0,1], \mu ) )$ (resp. $\{N_s: 0\leq s\leq 1 \}\subset B(L^2([0,1], \nu)))
 $ where $M_s$ (resp. $N_s$) is the projection onto the space $ L^2([0,s], \mu ) $
 (resp. $L^2([0,s], \nu ) $).

The map $\theta : \cl M_1\rightarrow \cl M_2: M_s\rightarrow N_{\gamma (s)}$ is 
a nest isomorphism so by Theorem \ref{2.11} the algebras $A=\Alg{\cl M_1}, B=\Alg{\cl M_2}$  
 are weakly$-*$ Morita equivalent. If the algebras  $A, B$ were $\Delta -$equivalent 
by Theorem \ref{3.7} there would exist unit vectors $\xi $ for $\cl M_1^{\prime\prime} ,$ 
  $\omega $ for $\cl M_2^{\prime\prime} $ and $\alpha \in AbsHom([0,1])$ such 
that $\alpha ( \phi _\xi (\cl M_1) )= \psi _\omega (\cl M_2) .$ 
From \cite[Proposition 7.22]{dav} we have that $m(\phi _\xi (\cl M_1) )=m(C)=0$ 
and since $m(\gamma (C))>0$ we have that $m(\psi _\omega (\cl M_2) )>0.$ This is a contradiction. $\qquad \Box$   

\section{A counterexample in Morita equivalence}

In this section we shall use the notions of TRO equivalence, of $\Delta -$equivalence, of 
stable isomorphism and we shall consider nest  and CSL algebras. See the appropriate definitions in sections 1, 3 and 4. 
If $C$ and $D$ 
are \textbf{unital} operator algebras which are strongly Morita equivalent then for 
every $\epsilon >0$ there exists a completely bounded isomorphism from $C\otimes _{min}\cl K$ 
onto $D\otimes _{min}\cl K$ with $\|\rho \|_{cb}<1+\epsilon $ and  
$\|\rho^{-1} \|_{cb}<1+\epsilon ,$ where $\cl K$ is the $C^*-$algebra of compact operators 
on a separable infinite dimensional Hilbert space $H$ and $\otimes _{min}$ is the spatial 
tensor product \cite[Corollary 7.10]{bmp}. It follows that for every $\epsilon >0$ 
there exists a completely bounded $w^*-$continuous isomorphism $\sigma $ from 
$C^{**} \stackrel{-}{\otimes} B(H)$ 
onto $D^{**} \stackrel{-}{\otimes}B(H)$  with $\|\sigma  \|_{cb}<1+\epsilon $ and  
$\|\sigma^{-1} \|_{cb}<1+\epsilon,$ where $\stackrel{-}{\otimes} $ is the normal spatial 
tensor product. One can wonder now, if the operator algebras $C^{**}$ and $D^{**}$ are 
stably isomorphic.

In this section we give a negative answer to this question. We present a counterexample of 
unital strongly Morita equivalent algebras $C$ and $D$ whose second duals are not 
stably isomorphic. Also for the algebras $C^{**}$ and $D^{**}$  there exist normal completely 
isometric representations  $\phi $ and $\psi $  respectively such that 
for every $\epsilon >0$ there exists an invertible bounded operator $T_\epsilon $ satisfying 
$\|T_\epsilon \|<1+\epsilon ,\;\; \|T_\epsilon ^{-1}\|<1+\epsilon ,\;\;$
$\phi (C^{**})=T_\epsilon ^{-1}\psi (D^{**})T_\epsilon $  and 
$\phi (C)=T_\epsilon ^{-1}\psi (D)T_\epsilon .$

Two nests $\cl N, \cl M$ acting on the separable Hilbert spaces $H, K$ respectively are 
called \textbf{similar} if there exists an order isomorphism $\theta : \cl N\rightarrow \cl M$ 
which preserves dimension of intervals. We say that an invertible operator $S\in B(H,K)$ 
implements $\theta $ if $\theta (N)$ is the projection onto $SN(H)$ for all $N\in \cl N.$ 
In what follows if $C$ is an operator algebra, $\Delta (C)$ is its diagonal $C\cap C^*.$

We fix similar nests $\cl N, \cl M$ as above with corresponding nest algebras 
$A=\mathrm{Alg}(\cl N)$ and $B=\mathrm{Alg}(\cl M)$ such that $\Delta (A)$ 
is a totally atomic maximal abelian selfadjoint algebra (\textbf{masa} in sequel) 
and  $\Delta (B)$ is a masa with a nontrivial continuous part, \cite[example 13.15]{dav}.
Suppose that $\theta : \cl N\rightarrow \cl M$ is an order isomorphism 
implementing similarity for $\cl N, \cl M.$  We denote by $A_0$ (resp. $B_0$) 
the algebra of compact operators belonging to $A$ (resp. $B$)
 and by $A_1$ (resp. $B_1$) the operator algebra $A_0+\mathbb{C}I_H$ (resp. 
$B_0+\mathbb{C}I_K$ ). We denote by $X$ the space $Op(\theta )$ and 
by $Y$ the space $Op(\theta^{-1} ).$ 

\begin{theorem}\label{5.1}\cite[Theorem 13.20]{dav}(Davidson) For every $\epsilon >0$ there exists an 
invertible bounded operator $S_\epsilon$ which implements $\theta $ such 
that $ \|S_\epsilon \|<1+\epsilon ,  \| S_\epsilon^{-1} \|<1+\epsilon .$ 
(Observe that $S_\epsilon \in X$ and $S_\epsilon^{-1} \in Y$ for all $\epsilon >0.$)
\end{theorem}

Suppose that $j: A_1\rightarrow A_1^{**}$ is the canonical embedding. We denote 
by $J_A$ the space $\overline{j(A_0)}^{w^*}.$

\begin{lemma}\label{5.2} 

(i) $A_1^{**}=J_A+\mathbb{C}I$

(ii)  $J_A\cap \mathbb{C}I=0.$

\end{lemma}
\textbf{Proof} 

(i) Since $|\lambda |\leq \|a+\lambda I_H\|$ for all compact operators $a$ the functional 
$$\rho : A_1\rightarrow \mathbb{C}:  a+\lambda I_H\rightarrow \lambda $$ belongs to 
$A_1^*.$ If $x\in A_1^{**}$ by the Goldstine Theorem there exists a net 
$(a_i+\lambda _iI_H)\subset A_0+\mathbb{C}I
$ converging in $w^*-$topology to $x.$ Since $(\lambda _i)$ converges to $\rho (x)$ 
we have that $(a_i)$ converges to $a\in J_A$ and so $x=a+\rho (x)\in J_A+\mathbb{C}I.$

(ii) Since $\rho |_{A_0}=0$ if $\lambda I\in J_A$ then $\lambda =0.$ 
So $J_A\cap \mathbb{C}I=0.\qquad \Box$

\bigskip

Suppose that $\iota : A_0\rightarrow A_0^{**}$ is the canonical embedding. In lemma 
\ref{2.9} we have proved that there exists a $w^*-$continuous completely isometric 
onto homomorphism  $\phi : A\rightarrow A_0^{**}$ extending $\iota .$ 

The map 
$\phi |_{A_1}: A_1\rightarrow A_0^{**}$ extends to a $w^*-$continuous completely 
contractive map $ \stackrel{\wedge }{\phi } : A_1^{**}\rightarrow A_0^{**}$ satisfying  
$ \stackrel{\wedge }{\phi } (j(a))=\phi (a)$ for all $a\in A_1.$ Also the 
 completely contractive map $j|_{A_0}: A_0\rightarrow A_1^{**}$ extends to 
a $w^*-$continuous completely 
contractive map $ \stackrel{\wedge }{\kappa  } : A_0^{**}\rightarrow A_1^{**} $ 
such that $ \stackrel{\wedge }{\kappa  } (\iota (a))=j(a)$ for all $a\in A_0.$ 
So the map $ \stackrel{\wedge }{\phi } \circ  \stackrel{\wedge }{\kappa  } : A_0^{**}\rightarrow 
A_0^{**} $ satisfies $$\stackrel{\wedge }{\phi }\circ  \stackrel{\wedge }{\kappa  } 
(\iota (a))=\stackrel{\wedge }{\phi } (j(a))=\phi (a)=\iota (a)$$ for all $a\in A_0.$    
It follows that $\stackrel{\wedge }{\phi }\circ \stackrel{\wedge }{\kappa }=id_{A_0^{**}}.$ 
Therefore $ \stackrel{\wedge }{\kappa } $ is a complete isometry.

We denote by $\theta $ the $w^*-$continuous completely isometric homomorphism  
$ \stackrel{\wedge }{\kappa } \circ 
\phi: A\rightarrow A_1^{**}.  $ Observe that
$$\theta (A)=\stackrel{\wedge }{\kappa } (\phi (A))=\stackrel{\wedge }{\kappa } (\overline{
\iota (A_0)}^{w^*})=\overline{j(A_0)}^{w^*}=J_A.$$
 Suppose that $p$ is the projection $\theta (id_A).$ Lemma \ref{5.2} implies that  $p^\bot \neq 0$ and 
$A_1^{**}= J_A\oplus \mathbb{C}p^\bot .$

\begin{lemma}\label{5.3} The algebra $A_1^{**}$ is completely isometric and $w^*-$
continuously isomorphic with the algebra $A\oplus \mathbb{C}$ acting on the Hilbert 
space $ H\oplus \mathbb{C} .$ 
\end{lemma}
\textbf{Proof} We define the map $\theta $ and the projection $p$ as in the 
above discussion. We define the completely isometric normal representation 
$$\pi : A_1^{**}=J_A\oplus \mathbb{C}p^\bot \rightarrow B(H\oplus \mathbb{C}): a\oplus \lambda 
p^\bot \rightarrow \theta ^{-1}(a)\oplus \lambda  $$ which is onto $A\oplus \mathbb{C}.
\qquad \Box$

\bigskip

For every $\epsilon >0$ we denote by $T_\epsilon $ the bounded 
invertible operator $S_\epsilon \oplus id_{ \mathbb{C} }\in B(H\oplus \mathbb{C} ,
K\oplus \mathbb{C} ).$ Also we denote the spaces $U=X\oplus \mathbb{C} \subset 
B(H\oplus \mathbb{C} ,K\oplus \mathbb{C} )$ and $V=Y\oplus \mathbb{C} \subset 
B(K\oplus \mathbb{C} ,H\oplus \mathbb{C} ).$ Observe that 
$U$ is a $B\oplus \mathbb{C} -A\oplus \mathbb{C} $ bimodule 
and $V$ is an $A\oplus \mathbb{C} -B\oplus \mathbb{C} $ bimodule.

By the above lemma $\pi (A_1^{**})=A\oplus \mathbb{C} .$ If $j: A_1\rightarrow A_1^{**}$ 
is the canonical embedding we have $\pi (j(a))=a\oplus 0$ for all 
$a\in A_0$ and $\pi (j(id_{A_1}))= id_{H\oplus 
\mathbb{C} }.$ So $$\pi (j(A_1))=\mathrm{span}\{a\oplus 0, id_{H\oplus 
\mathbb{C} }, a\in A_0\}.$$

Similarly if $j_2: B_1\rightarrow B_1^{**}$ is the canonical embedding  there exists a normal completely 
isometric onto homomorphism $\rho : B_1^{**}\rightarrow B\oplus \mathbb{C} $ 
such that $$\rho (j_2(B_1))=\mathrm{span}\{b\oplus 0, id_{K\oplus 
\mathbb{C} }, b\in B_0\}.$$

Since $S_\epsilon ^{-1}B_0S_\epsilon =A_0$ and $S_\epsilon ^{-1}BS_\epsilon =A$ 
we have that $$ T_\epsilon ^{-1}\rho (j_2(B_1))T
_\epsilon =\pi (j(A_1)),\;\; T_\epsilon ^{-1}\rho (B_1^{**})T
_\epsilon = \pi (A_1^{**})$$ for all $\epsilon >0.$ 

In the following lemmas \ref{mor}, \ref{5.4} we identify the algebra $A_1^{**}$ with $A\oplus \mathbb{C}, $ 
the algebra $B_1^{**}$ with $B\oplus \mathbb{C}, $ the algebra 
$A_1$ with $\pi (j(A_1))$ and the algebra 
$B_1$ with $\rho (j_2(B_1)). $

\begin{lemma}\label{mor} The algebras $ A_1^{**}$ and $B_1^{**} $ are weakly$-*$ Morita equivalent. 
\end{lemma}
\textbf{Proof} Let $U, V$ and $T_\epsilon, \epsilon >0 $ be as in the above discussion. 
The completely contractive bilinear map $V\times U\rightarrow A_1^{**} : (v,u)\rightarrow vu$  
is separately $w^*$-continuous, $ B_1^{**} -$balanced and $A_1^{**} -$module map. 
So induces the $w^*$-continuous completely contractive and $A_1^{**} -$module map
 $$\tau : V\otimes ^{\sigma h}_{ B_1^{**} }U \rightarrow 
A_1^{**} : v \otimes _{B_1^{**} } u\rightarrow vu.$$
We shall prove that $\tau $ is isometric: If $(v_i)\subset V, (u_i)\subset U$ and $\epsilon >0$ we have:
$$ \nor{\sum_{i=1}^nv_i\otimes _{B_1^{**} } u_i} = 
\nor{\sum_{i=1}^n(T_\epsilon ^{-1}T_\epsilon v_i)\otimes _{B_1^{**} } u_i} . $$
Since $T_\epsilon v_i\in UV\subset B_1^{**} $ the last norm is equal with 
\begin{align*}&\nor{ \sum_{i=1}^n  T_\epsilon ^{-1} \otimes _{B_1^{**} } (T_\epsilon v_iu_i)}=
\nor{( T_\epsilon ^{-1} \otimes _{B_1^{**} } T_\epsilon )( \sum_{i=1}^n v_iu_i )}\\ 
&\leq \|T_\epsilon ^{-1} \|\|T_\epsilon \| \nor{\sum_{i=1}^n v_iu_i } \leq (1+\epsilon )^2
 \nor{\sum_{i=1}^n v_iu_i } . 
\end{align*}
We let $\epsilon \rightarrow 0$ and we have that 
$$ \nor{\sum_{i=1}^nv_i\otimes _{ B_1^{**} } u_i} =\nor{\sum_{i=1}^n v_iu_i } .$$
Similarly we can prove that $\tau $ is completely isometric. Since $A=\wsp{YX}$ 
we have that $ A_1^{**} =\wsp{VU}$ and so by the Krein-Smulian Theorem $\tau $ is onto $A_1^{**} .$ 
The proof of the fact $B_1^{**} \cong U\otimes ^{\sigma h}_{ A_1^{**} }V $  is similar.
$\qquad \Box$

\begin{lemma}\label{5.4} The algebras $A_1$ and $B_1$ are strongly Morita equivalent. 
\end{lemma}  
\textbf{Proof} It suffices to prove that they have equivalent 
categories of left operator modules \cite{blecher}. If $C$ is an operator 
algebra we denote by $ \;_C mod $ the category of left operator modules over $C.$ 
We assume that every $Z\in \;_C mod $ is essential, i.e. the linear 
span of $CZ$ is dense in $Z.$ If  $Z_1, Z_2\in \;_C mod $ the space 
of morphisms $Hom_C(Z_1, Z_2)$ is the space of completely bounded maps 
$F: Z_1\rightarrow Z_2$ which are $C-$module maps.

We fix an operator $T=T_{\epsilon _0}$ for 
$\epsilon _0>0.$ If $Z\in \;_{A_1}mod$ then $Z^{**}$ is a left dual operator module over $A_1^{**}$ in a canonical 
way \cite[3.8.9]{bm}. We denote by $\cl F(Z)$ the subspace of $U\otimes ^{\sigma h}_{A_1^{**}}
Z^{**}$ 
$$\cl F(Z)=\lsp{ Ta\otimes _{A_1^{**}}z : a\in A_1, z\in Z}.$$
Since $U\otimes ^{\sigma h}_{A_1^{**}}
Z^{**}$ is a left operator module over $B_1^{**}$ and 
$$b(Ta\otimes _{A_1^{**}}z)= (bTa)\otimes _{A_1^{**}}z = T( T^{-1}bT a)\otimes 
_{A_1^{**}}z $$ with $T^{-1}bT\in A_1$ for all $b\in B_1,$ $\cl F(Z)$ is a left 
operator $B_1-$module. 

If $W\in \;_{B_1}mod$ we denote by $\cl G(W)$ the subspace of $V\otimes ^{\sigma h}_{B_1^{**}}
W^{**}$ 
$$\cl G(W)=\lsp{ aT^{-1}\otimes _{B_1^{**}}w : a\in A_1, w\in W}.$$
Since $V\otimes ^{\sigma h}_{B_1^{**}}
W^{**}$ is a left operator module over $A_1^{**},$ 
clearly $\cl G(W)\in \;_{\cl A_1}mod.$ 

Now $$ \cl G(\cl F(Z)) =\lsp{ a_2T^{-1} \otimes _{B_1^{**}} Ta_1 \otimes _{A_1^{**}} z
 : a_1, a_2\in A_1, z\in Z}$$ is a left operator module over $A_1$ 
and subspace of the space $ V \otimes ^{\sigma h}_{B_1^{**}} U\otimes ^{\sigma h}_{A_1^{**}} 
Z^{**} .$ The $w^*-$Morita equivalence $A_1^{**} \cong V\otimes ^{\sigma h}_{ B_1^{**} }U,
 $ $B_1^{**} \cong U\otimes ^{\sigma h}_{ A_1^{**} }V $ induces (\cite[Theorem 3.5]{bk}) 
 a complete isometry  
$$V \otimes ^{\sigma h} _{B_1^{**}} U\otimes ^{\sigma h}_{A_1^{**}} 
Z^{**} \rightarrow Z^{**}: v\otimes_{B_1^{**}}  u\otimes_{A_1^{**}}  z\rightarrow vuz$$
which restricts to a  completely isometric map 
$$R_Z: \cl G(\cl F(Z)) \rightarrow Z: a_2T^{-1} \otimes _{B_1^{**}} Ta_1 \otimes _{A_1^{**}} z
\rightarrow a_2a_1z$$ 
for all $a_1, a_2\in A_1, z\in Z.$ This map is clearly onto $Z.$

Every morphism $F\in Hom_{A_1}(Z_1, Z_2)$ can be extended to a morphism 
$\stackrel{\wedge }{F}$ belonging to  $Hom^\sigma _{A_1^{**}}(Z_1^{**}, Z_2^{**}),$ 
the space of $w^*-$continuous completely bounded $A_1^{**}-$module maps. (Use for example 
\cite[1.4.8]{bm}).

 The weak$-*$ Morita equivalence  
$A_1^{**}\cong V\otimes ^{\sigma h}_{B_1^{**}}U,\;\;  
B_1^{**}\cong U\otimes ^{\sigma h}_{ A_1^{**} }V,$ generates (\cite[Theorem 3.5]{bk}) a normal completely contractive functor 
$\stackrel{\wedge }{\cl F} $ between the left dual operator modules of $A_1^{**} $ and $B_1^{**} $ such that 
$$ \stackrel{\wedge }{\cl F} ( \stackrel{\wedge }{F}  ): U\otimes ^{\sigma h}_{ A_1^{**} }Z_1^{**}
\rightarrow U\otimes ^{\sigma h}_{A_1^{**}}Z_2^{**}: u\otimes _{A_1^{**} }z \rightarrow 
u\otimes _{A_1^{**} }\stackrel{\wedge }{F}  (z). $$
Since $$  \stackrel{\wedge }{\cl F} ( \stackrel{\wedge }{F}  )( Ta\otimes _{A_1^{**}}z )=Ta\otimes _{A_1^{**}}F(z) $$ for all 
$a\in A_1, z\in Z_1$ the operator 
 $\stackrel{\wedge }{\cl F} ( \stackrel{\wedge }{F})$ maps $\cl F(Z_1)$ into $\cl F(Z_2).$ 
So we can define $$ \cl F(F) = \stackrel{\wedge }{\cl F} ( \stackrel{\wedge }{F}) |_{ \cl F(Z_1)  
}:  \cl F(Z_1) \rightarrow \cl F(Z_2). $$
 We can easily check that $\cl F(F) 
\in Hom_{B_1}(\cl F(Z_1), \cl F(Z_2)). $

In this way we define functors $\cl F: \;_{A_1}mod\rightarrow \;_{B_1}mod$
  and  $\cl G: \;_{B_1}mod\rightarrow \;_{A_1}mod.$ Using the above complete isometries $\{R_Z: Z\in \;_{A_1}mod \}$ 
we can prove that the functor $\cl G\cl F$ is equivalent to the identity functor 
$1_{\;_{A_1}mod }$ and the functor $\cl F\cl G$ is equivalent to the identity functor 
$1_{\;_{B_1}mod }. \qquad \Box$

\begin{theorem}\label{5.5} Strong Morita equivalence of unital operator algebras 
doesn't imply $\Delta -$equivalence of the second dual operator algebras.
\end{theorem}
\textbf{Proof} We recall the unital operator algebras $A_1, B_1$ which are strongly Morita equivalent 
by the above  lemma. We shall prove that the algebras $A_1^{**}, B_1^{**}$ are not 
$\Delta -$equivalent. Suppose that they are $\Delta -$equivalent. We define the completely 
isometric normal representation (see Lemma \ref{5.3}) 
$$\pi : A_1^{**}\rightarrow B(H\oplus \mathbb{C}): a\oplus \lambda p^\bot \rightarrow \theta 
^{-1}(a)\oplus \lambda .$$

The algebra $\pi (A_1^{**}) =A\oplus \mathbb{C}$ is a CSL algebra with lattice 
$$\{N\oplus 0, N\oplus \mathbb{C}: N\in \cl N\}.$$
Suppose that $B_1^{**}=J_B\oplus \mathbb{C}q^\bot $ where $q$ is the identity 
of the algebra $J_B$ and $J_B$ is isomorphic with the algebra $B.$ By \cite[Theorem 2.7]{ele2} 
there exists a completely isometric normal representation $\sigma $ of $B_1^{**}$ on a Hilbert space 
 $K_1\oplus K_2$ of the form $\sigma (b\oplus \lambda q^\bot )=\sigma _1(b)\oplus 
\lambda I_{K_2}$ for all $b\in J_B, \lambda \in \mathbb{C}$ such that the algebras $\pi (A_1^{**}), \sigma (B_1^{**})  $
 are TRO equivalent. Since  $ \pi (A_1^{**}) $ is a CSL algebra, $\sigma (B_1^{**})  $ 
is also a CSL algebra, \cite[Remark 5.5]{ele}. So the algebra $ \sigma (B_1^{**})  $ contains 
a masa. It follows that $dim K_2=1.$ So we may assume that $\sigma (B_1^{**})  $ 
is a CSL algebra acting on $K_1\oplus \mathbb{C}.$ 

Since $\Delta (A)$ (resp. $\Delta (B)$) is a masa, then $\Delta (\pi (A_1^{**}) )$ 
(resp. $\Delta ( \sigma (B_1^{**})  )$) is also a masa. The algebras 
$\Delta (\pi (A_1^{**}) ), \Delta ( \sigma (B_1^{**})  )$ are TRO equivalent 
\cite[Proposition 2.5]{ele}.
 But TRO equivalence between masas is a unitary equivalence (use for example \cite[Theorem 3.2]
{ele}). 
This is a contradiction because 
$\Delta (\pi ( A_1^{**} ) )= \Delta (A)\oplus \mathbb{C} $ is a totally atomic masa 
and  the masa $\Delta ( \sigma (B_1^{**})  )\cong  \Delta (B)\oplus \mathbb{C} $  has 
a nontrivial continuous part. So the algebras $A_1^{**}, B_1^{**} $ are not 
$\Delta -$equivalent.$\qquad \Box$

\bigskip

{\em Acknowledgement:} I wish to thank Prof.~D.~Blecher who pointed out the 
problem in section 5 to me.


\begin{thebibliography}{99}

\bibitem{arv}
W. B. Arveson, {\em Operator algebras and invariant subspaces}, Ann. of Math.
100 (1974) 433-532.

\bibitem{blecher}
D. P. Blecher, {\em A Morita theorem for algebras of operators on Hilbert space},  
J. Pure Appl. Algebra 156 (2001), 156-169.



\bibitem{bk}
D. P. Blecher and U. Kashyap, {\em Morita equivalence of dual
operator algebras}, J. Pure Appl. Algebra, 212 (2008), 2401-2412.

\bibitem{bkraus}
D. P. Blecher and J. E. Kraus, {\em On a generalization of $W^*-$modules}, 
Arxiv: 0910.5404 

\bibitem{bm}
D. P. Blecher and C. Le Merdy, {\em Operator algebras and their
modules---an operator space approach}, Oxford University Press,
2004

\bibitem{bmp} D. P. Blecher, P. S. Muhly and V. I. Paulsen, {\em Categories of Operator Modules
(Morita Equivalence and Projective Modules)}, Mem. Amer. Math.
Soc. 143 (2000), no.681.


\bibitem{dav}
K. R. Davidson,
\newblock {\em Nest algebras},
\newblock Longman Scientific \& Technical, Harlow, 1988

\bibitem{dav2}
K. R. Davidson, {\em Commutative subspace lattices}, Indiana Univ. Math. J. 
 27 (3) (1978) 479-490.

\bibitem{er} E. Effros and Z.-J. Ruan, {\em Operator Spaces}, Clarendon Press, Oxford, 2000

\bibitem{ele}
G. K. Eleftherakis, {\em TRO equivalent algebras}, Houston J. of Mathematics, to appear.

\bibitem{ele1}
G. K. Eleftherakis, {\em A Morita type equivalence for dual
operator algebras}, J. Pure Appl. Algebra 212 (2008) no 5,
1060-1071

\bibitem{ele2}
G. K. Eleftherakis, {\em Morita type equivalences and reflexive
algebras}, J. Operator Theory, to appear

\bibitem{elepaul}
G. K. Eleftherakis and V. I. Paulsen, {\em Stably isomorphic dual
operator algebras}, Math. Ann. 341 (2008), no 1, 99-112

\bibitem{kashyap} U. Kashyap, A Morita theorem for dual operator
algebras, J. Funct. Analysis, 256 (2009) 3545-3567.

\bibitem{paul}
V. I. Paulsen, {\em Completely Bounded Maps and Operator
Algebras}, Cambridge University Press, 2002

\bibitem{pisier}
G. Pisier, {\em Introduction to Operator Space Theory}, Cambridge
University Press, 2003

\bibitem{rief} M. Rieffel, {\em Morita equivalence for C*-algebras and W*-algebras},
J. Pure Appl. Algebra 5 (1974), 51-96



\end{thebibliography}
\end{document}